\documentclass{amsart}

\usepackage{amssymb}
\usepackage[all]{xy}

\newtheorem{thm}{Theorem}

\newtheorem{cor}[thm]{Corollary}

\newtheorem{prop}[thm]{Proposition}

\newtheorem{conj}[thm]{Conjecture}
   
\theoremstyle{definition}

\newtheorem{say}[thm]{}
\newtheorem{exmp}[thm]{Example}

\newtheorem*{ack}{Acknowledgments}      

\newtheorem{defn-thm}[thm]{Definition--Theorem}  
\newtheorem{defn-lem}[thm]{Definition--Lemma}  

\newtheorem{reminder}[thm]{Reminder}

\theoremstyle{remark}

\renewcommand{\o}[0]{{\mathcal O}} 
\newcommand{\z}[0]{{\mathbb Z}}

\renewcommand{\r}[0]{{\mathbb R}} 

\renewcommand{\a}[0]{{\mathbb A}}

\newcommand{\q}[0]{{\mathbb Q}}

\newcommand{\qtq}[1]{\quad\mbox{#1}\quad}

\newcommand{\supp}[0]{\operatorname{Supp}}    
    
\newcommand{\codim}[0]{\operatorname{codim}}

\newcommand{\diff}[0]{\operatorname{Diff}}

\newcommand{\res}[0]{\operatorname{\mathcal R}}

\newcommand{\cl}[0]{\operatorname{Cl}}

\newcommand{\rup}[1]{\lceil{#1}\rceil}
\newcommand{\rdown}[1]{\lfloor{#1}\rfloor}

\newcommand{\onto}[0]{\twoheadrightarrow}

\newcommand{\simq}[0]{\sim_{\q}}

\newcommand{\simr}[0]{\sim_{\r}}

\newcommand{\depth}[0]{\operatorname{depth}} 
\newcommand{\tsum}[0]{\textstyle{\sum}}

\newcommand{\coeff}[0]{\operatorname{coeff}}

\def\loccoh#1.#2.#3.#4.{H^{#1}_{#2}(#3,#4)}

\DeclareMathAlphabet{\mathchanc}{OT1}{pzc}%
                                {m}{it}

\usepackage[all]{xy}\xyoption{dvips}

\newcommand{\congt}[0]{\operatorname{\cong_{tor}}}

\begin{document}
\bibliographystyle{amsalpha}

 \title{Log-plurigenera in stable families of surfaces}
 \author{J\'anos Koll\'ar}

\begin{abstract} We study the flatness of log-pluricanonical sheaves on stable families of surfaces.
\end{abstract}

 \maketitle

The paper  \cite{k-lpg1}  studies flatness of  the pluricanonical sheaves
$ \omega_{X/S}^{[m]}\bigl(\rdown{m\Delta}\bigr) $ 
for locally stable morphisms over reduced base schemes. Positive results are obtained for families with  normal  generic fibers, provided  every divisor appears in $\Delta$ with coefficient $\geq \frac12$.
While examples show that the bound $\geq \frac12$ is sharp, 
presumably the normality condition is not necessary. The aim of this note is to prove this for families of surfaces. 
Unfortunately, the proof relies on the classification of
slc surface pairs, thus it is unlikely to generalize to higher dimensions. 

\begin{thm}  \label{half.pgsh.2d.inv.thm}
Let $S$ be a reduced scheme over a field of characteristic 0 and $f:(X,\Delta) \to S$ a locally stable morphism of relative dimension 2. Assume  that $\coeff\Delta\subset [\frac12, 1]$.
Then, for every $m\in\z$  and $B\subset \rdown{\Delta}$, the sheaves
$$
\omega_{X/S}^{[m]}\bigl(\rdown{m\Delta}-B\bigr)
\eqno{(\ref{half.pgsh.2d.inv.thm}.1)}
$$
are flat over $S$ and commute with base change.
\end{thm}

{\it Warning \ref{half.pgsh.2d.inv.thm}.2.}  The Theorem and its Corollary  hold for every variant of local stability I know of if either $\coeff\Delta\subset (\frac12, 1]$ or if $S$ is unibranch. In general we need to assume also that 
$\rdown{\Delta_s}=\rdown{\Delta}_s$ holds for every $s\in S$. See
\cite[Sec.6]{k-lpg1} for  relevant examples and 
\cite[Chap.4]{k-modbook} for a detailed discussion of the issues. 

\medskip
As a first consequence we obtain that, if $S$ is connected, then  
 the Hilbert function of the fibers
$$
\chi\bigl(X_s, \omega_{X_s}^{[m]}\bigl(\rdown{m\Delta_s}\bigr)\bigr)
\eqno{(\ref{half.pgsh.2d.inv.thm}.3)}
$$
is  independent of  $s\in S$.  If  $f:(X,\Delta) \to S$ a  stable, that is, if
$K_{X/S}+\Delta$ is also $f$-ample, then  by Serre vanishing, the log plurigenera
$$
p_m(X_s, \Delta_s):=h^0\bigl(X_s, \omega_{X_s}^{[m]}\bigl(\rdown{m\Delta_s}\bigr)\bigr)
\eqno{(\ref{half.pgsh.2d.inv.thm}.4)}
$$
are also  independent of  $s\in S$ for $m\gg 1$. 
We can be more precise  if we restrict the coefficients further.

\begin{cor} \label{st.pg.inv.cor}
 Let $S$ be a reduced scheme over a field of characteristic 0 and $f:(X,\Delta) \to S$ a stable morphism of relative dimension 2  such that $\coeff\Delta\subset \{\frac12, \frac23, \frac34, \dots, 1\}$. Then, for every $m\geq 2$,
\begin{enumerate}
\item  $R^if_*\omega_{X/S}^{[m]}\bigl(\rdown{m\Delta}\bigr)=0$ for $i>0$ and
\item $f_*\omega_{X/S}^{[m]}\bigl(\rdown{m\Delta}\bigr)$ is locally free and commutes with base change.
\end{enumerate}
\end{cor}

Both the Therem and the Corollary  should hold in higher dimensions as well, hence the surface case is rather special. Therefore the main interest of this note may be the observation that the gluing theory of log pluricanonical sheaves  on slc pairs seems much more complicated than the   gluing  of  slc pairs themselves. The latter  was introduced in  \cite{k-source} and  discussed in detail in \cite[Chap.5]{kk-singbook}.
\medskip

In Section~\ref{sec.1} we reduce the Theorem to  a claim about
slc threefolds, which is then proved in Section~\ref{sec.2}. 
The proof uses detailed information about certain non-normal slc surfaces. 
These include a partial  classification of  non-normal slc surfaces, given  in   Section~\ref{sec.3}, and the  computation of  
the Poincar\'e residue map on their irreducible components, treated in   Section~\ref{sec.4}. Corollary~\ref{st.pg.inv.cor} is proved in 
Section~\ref{sec.5}.

\begin{ack} I thank   Chenyang~Xu for insightful comments. 
Partial  financial support    was provided  by  the NSF under grant number
 DMS-1362960.
\end{ack}

\section{Non-normal slc threefolds}\label{sec.1}

Using  \cite[Prop.16]{k-lpg1},  it is sufficient to prove Theorem~\ref{half.pgsh.2d.inv.thm}  when
$S$ is regular and of dimension 1. Furthermore, the latter
is equivalent to proving the following variant, which will be the focus of our attention from now on.

\begin{prop} \label{tech.prop.2d}
Let $(x\in X,H+\Delta)$ be a 3-dimensional, local   slc pair  over a field of characteristic 0 where $H$ is Cartier and $\coeff\Delta\subset [\frac12, 1]$. 
Then  $$
\depth_x\omega_X^{[m]}\bigl(\rdown{m\Delta}-B\bigr)=3
\qtq{for every $m\in \z$ and $B\subset \rdown{\Delta}$.}
$$
Equivalently, $\omega_X^{[m]}\bigl(\rdown{m\Delta}-B\bigr)$ satisfies Serre's condition $S_3$. 
\end{prop}

For a coherent sheaf whose support has  dimension $\leq 3$, being $S_3$ is equivalent to being Cohen-Macaulay. If $\dim X\geq 4$ then the sheaves
$\omega_X^{[m]}\bigl(\rdown{m\Delta}-B\bigr)$ are frequently not Cohen-Macaulay, but a (slight modification of) the $S_3$ condition is expected to hold; see  \cite[Prop.5]{k-lpg1}. This is why we state Proposition~\ref{tech.prop.2d} using the $S_3$ condition. 

The method of  \cite{k-lpg1}, which proves  Proposition~\ref{tech.prop.2d} in case $X$ is normal,  has 3 steps. The first, going back to 
\cite{ale-lim, k-dep} establishes the case when
$mK_{X}+\rdown{m\Delta}-B$ is $\q$-Cartier. The second constructs a
small modification $\pi:X'\to X $ such that
$mK_{X'}+\rdown{m\Delta'}-B'$ is $\q$-Cartier and the third uses $X'$ to obtain the conclusion.  As observed in  \cite[Exmp.22]{k-lpg1}, the second step usually does not hold if $X$ is not normal; there are obstructions in codimension 2 and also in higher codimensions.
In this note  we deal with the codimension 2 obstruction. In the theory of slc pairs,  the higher codimension  obstructions usually  behave quite diferently, and there are several instances when the higher codimension case is easier. So there is some reason to believe that handling the  codimension 2 obstruction may be a useful step in general.

 If $X$ is normal  then the conclusion of Proposition~\ref{tech.prop.2d} is known to hold in all dimensions by
\cite[Prop.5]{k-lpg1}. 
   Thus it remains to understand what happens when $X$ is
non-normal.  The gluing method of \cite{k-source} suggests that
one should be able to treat $X$ by first working on its normalization $\bigl(\bar X, \bar D+ \bar\Delta\bigr)$, then proving 
 compatibility with the gluing involution $\tau$ 
 and finally  desceding to $X$; see \cite[Chap.5]{kk-singbook} for details.
Compatibility  with the gluing involution turns out to be quite subtle.
There are  2 variants: 
\begin{itemize}
\item divisor version, working with $mK_{\bar X}+m\bar D+\rdown{m\bar\Delta}$ and the different, and 
\item sheaf version,   working with $\omega_{\bar X}^{[m]}\bigl(m\bar D+\rdown{m\bar\Delta}\bigr)$ and the Poincar\'e residue maps
$\res_{\bar X/\bar D}^m$ as in \cite[Sec.4.1]{kk-singbook}.
\end{itemize}
Unexpectedly, the 2 variants  are {\em not} equivalent, and
 compatibility fails for both of them. For the  divisor version see
Example~\ref{pre.key.plt.exmp}, for the sheaf version see
Examples ~\ref{diff.nonr.exmp.1} and \ref{diff.nonr.exmp.2}.
However,  the  sheaf version does hold in many instances
and one can describe quite well all cases when it fails.

With this in mind, 
first we focus on $H$ and prove a rather complete
 \'etale-local classification of such non-normal surface pairs
$\bigl(H, \diff_H\Delta\bigr)$ in
Theorem~\ref{slc.cent.12.prop}. 
This in turn implies the following
description of the pair $(X,H+\Delta)$.
This classification also shows that $x\in \supp B$ can happen only in the simpler case (\ref{codim.2.lccent.prop}.1), thus we can mostly ignore $B$ in the sequel.

\begin{prop} \label{codim.2.lccent.prop}
Let $(x\in X,H+\Delta)$ be a 3-dimensional, strictly Henselian,  slc pair  over a field of characteristic 0 where $H$ is Cartier and $\coeff\Delta\subset [\frac12, 1]$. 
Assume that  $X$ is not normal. 
Then  one of the following holds.
\begin{enumerate}
\item The point $x$ is an lc center and $2(K_X+H+\Delta)$ is Cartier at  $x$.
\item  The point $x$ is not an lc center  and  $X$   has 2 irreducible components
$(x_i\in X_i, D_i+H_i+\Delta_i)$ where $D_i$ denotes the conductor. Furthermore, the $D_i$ are smooth and 
the Poincar\'e residue maps
$$
\res_{X_i/D_i}^m: 
\omega_{X_i}^{[m]}\bigl(mD_i+mH_i+\rdown{m\Delta_i}\bigr)\to
\omega_{D_i}^{[m]}\bigl(\rdown{m\diff_{D_i}(H_i+\Delta)}\bigr)
$$
are surjective for every $m$.
\item  The point $x$ is not an lc center, $X$ is irreducible  and  it has a quasi-\'etale double cover as in {\rm (2)}. 
\end{enumerate}
\end{prop}

\section{Proof of the main results}\label{sec.2}

\begin{say}[A reformulation  of Proposition~\ref{tech.prop.2d}]
\label{tech.prop.2d.v2}
The Poincar\'e residue map
$$
\res^m_{X/H}: \omega_X^{[m]}\bigl(mH+\rdown{m\Delta}-B\bigr)\to
\omega_H^{[m]}\bigl(\rdown{m\diff_H\Delta}-B|_H\bigr)
\eqno{(\ref{tech.prop.2d.v2}.1)}
$$
can be factored through the injection
$$
\omega_X^{[m]}\bigl(mH+\rdown{m\Delta}-B\bigr)|_H\to
\omega_H^{[m]}\bigl(\rdown{m\diff_H\Delta}-B|_H\bigr),
\eqno{(\ref{tech.prop.2d.v2}.2)}
$$
which is an isomorphism on $H\setminus\{x\}$, where both sheaves are locally free. Thus we see that 
$$
\depth_x\omega_X^{[m]}\bigl(mH+\rdown{m\Delta}-B\bigr)=3
\quad\Leftrightarrow\quad \res^m_{X/H}\qtq{is surjective.}
\eqno{(\ref{tech.prop.2d.v2}.3)}
$$
\end{say}

\begin{say}[Proof of Proposition~\ref{codim.2.lccent.prop}] \label{codim.2.lccent.prop.pf}
It is easy to establish that   $\coeff\bigl(\diff_H\Delta\bigr)\subset [\frac12, 1]$, see for instance  \cite[3.45]{kk-singbook}.

If $x$ is an lc center of $(X,H+\Delta)$  then it is also 
 an lc center of $\bigl(H, \diff_H\Delta\bigr)$ by adjunction \cite[4.9]{kk-singbook}. Thus Theorem~\ref{slc.cent.12.prop} shows that
$2(K_H+\diff_H\Delta)$ is Cartier.  Therefore  $ 2(K_X+\Delta)$ is also  Cartier 
by \cite[XIII]{sga2} or \cite[2.90]{k-modbook}, giving (\ref{codim.2.lccent.prop}.1).

If $x$ is not an lc center of $(X,H+\Delta)$  then it is also 
not  an lc center of $\bigl(H, \diff_H\Delta\bigr)$ by adjunction.
Thus $\bigl(H, \diff_H\Delta\bigr)$ is as described in (\ref{slc.cent.12.prop}.2). In particular, $X$ has 1 or 2 irreducible components. If $X$ has only 1 irreducible component then 
by  \cite[5.23]{kk-singbook}  it has a quasi-\'etale double cover
with 2 irreducible components. This gives case (\ref{codim.2.lccent.prop}.3).

It remains to consider the case when  $X$ has  2 irreducible components.
Then $C_i:=D_i\cap H_i$ is smooth by (\ref{slc.cent.12.prop}.2) hence the $D_i$ are smooth. The various residue maps sit in a diagram 
$$
\begin{array}{ccc}
\omega_{X_i}^{[m]}\bigl(mD_i+mH_i+\rdown{m\Delta_i}\bigr)
& \stackrel{\res_{X_i/D_i}^m}{\longrightarrow} &
\omega_{D_i}^{[m]}\bigl(\rdown{m\diff_{D_i}(H_i+\Delta)}\bigr)\\[1ex]
\res_{X_i/H_i}^m\downarrow\hphantom{\res_{X_i/H_i}^m} && \hphantom{\res_{H_i/C_i}^m}\downarrow\res_{H_i/C_i}^m\\
\omega_{H_i}^{[m]}\bigl(\rdown{m(D_i+\Delta_i)}\bigr)
& \stackrel{\res_{H_i/C_i}^m}{\longrightarrow} &
\omega_{C_i}^{[m]}\bigl(\rdown{m\diff_{C_i}(H_i+\Delta)}\bigr)
\end{array}
$$
(Note that, as discussed in \cite[4.18]{kk-singbook},  this diagram commutes for even $m$ but only commutes up to sign if $m$ is odd. This has no bearing on which maps are  surjective.)
Since $X_i$ and $D_i$ are normal, the sheaves
$$
\omega_{X_i}^{[m]}\bigl(mD_i+mH_i+\rdown{m\Delta_i}\bigr) \qtq{and} \omega_{D_i}^{[m]}\bigl(\rdown{m\diff_{D_i}(H_i+\Delta)}\bigr)
$$
are $S_3$  by \cite[Prop.5]{k-lpg1}, thus the vertical arrows are surjective by (\ref{tech.prop.2d.v2}.3). 
Furthermore,  (\ref{key.plt.exmp}.1) shows that the bottom horizontal arrow
$\res_{H_i/C_i}^m $ is surjective. Thus the Nakayama lemma implies that 
the top horizontal arrow
$\res_{X_i/D_i}^m $ is also surjective, proving (\ref{codim.2.lccent.prop}.2).
 \qed
\end{say}

\begin{say}[Proof of Proposition~\ref{tech.prop.2d}]
The case when $X$ is normal at $x$ is proved in 
\cite[Prop.5]{k-lpg1},
thus assume form now on that $X$ is not normal at $x$.
 Proposition~\ref{codim.2.lccent.prop} classifies these into 3 cases, we consider them separately.

   In case (\ref{codim.2.lccent.prop}.1) 
$\coeff\Delta\subset \{\frac12, 1\}$ and 
 $\omega_X^{[2]}\bigl(\rdown{2\Delta}\bigr)= \omega_X^{[2]}\bigl({2\Delta}\bigr)$ is free. 
So, if  $m$ is even then
$\omega_X^{[m]}\bigl(\rdown{m\Delta}-B\bigr)\cong \o_X(-B)$ and if
$m$ is odd then 
$\omega_X^{[m]}\bigl(\rdown{m\Delta}-B\bigr)\cong \omega_X\bigl(\rdown{\Delta}-B\bigr)$.  We can now use 
\cite[7.20]{kk-singbook} first for $B\leq \Delta$ and then for
$-K_X-\rdown{\Delta}+B\simq \{\Delta\}+B\leq \Delta$ to conclude that
$\depth_x\omega_X^{[m]}\bigl(\rdown{m\Delta}-B\bigr)=3$. 
(Note that we could also have used the more general 
Corollary~\ref{tech.prop.std.cor}. Both of these arguments  apply in all dimensions.) 

 In case (\ref{codim.2.lccent.prop}.2)  the $H_i$ do not have any further role, so set $\Theta_i:=H_i+\Delta_i$. 
We know that  $X$ is obtained by gluing the 2 components
$(X_i, D_i+\Theta_i)$  using an  isomorphism
$\tau: \bigl(D_1, \diff_{D_1}\Theta_1\bigr)\cong  \bigl(D_2, \diff_{D_2}\Theta_2\bigr)$. Thus, by
\cite[5.8]{kk-singbook},  a pair of  sections 
$\sigma_i$ of $ \omega_{X_i}^{[m]}\bigl(\rdown{m\Theta_i}+mD_i\bigr)$
glues to a section of
$\omega_X^{[m]}\bigl(\rdown{m\Theta}\bigr)$ iff
$$
\res_{X_1/D_1}(\sigma_1)= (-1)^m\cdot \tau^*\bigl(\res_{X_2/D_2}(\sigma_2)\bigr).
$$
Equivalently,  we have an exact sequence
$$
0\to \omega_X^{[m]}\bigl(\rdown{m\Theta}\bigr)
\to 
\bigoplus_{i=1,2}\omega_{X_i}^{[m]}\bigl(\rdown{m\Theta_i}+mD_i\bigr)
\stackrel{R_D}{\longrightarrow}
\omega_{D_1}^{[m]}\bigl(\rdown{m\diff_{D_1}\Theta_1}\bigr)
\to 0,
$$
where $R_D=\res_{X_1/D_1}-\ (-1)^m\cdot \tau^*\circ \res_{X_2/D_2}$. 
The sheaves in the middle and  on the right are   $S_3$  by  \cite[Prop.5]{k-lpg1}. Thus the 
sheaf on the left is  also $S_3$; cf.\ \cite[2.60]{kk-singbook}.

In case (\ref{codim.2.lccent.prop}.3) let
$\pi:  (\tilde X, \tilde H+\tilde\Delta)\to (X, H+\Delta)$ be the double cover.
We already proved that 
 $\omega_{\tilde X}^{[m]}\bigl(\rdown{m\tilde\Delta}-\tilde D\bigr)$ is $S_3$,
hence so is  $\omega_X^{[m]}\bigl(\rdown{m\Delta}-D\bigr)$,  which is a direct summand of $\pi_*\omega_{\tilde X}^{[m]}\bigl(\rdown{m\tilde\Delta}-\tilde D\bigr)$. \qed
\end{say}

\section{Non-normal surface pairs}\label{sec.3}

Next we describe the pairs $(H, \diff_H\Delta)$ that arise in 
Proposition~\ref{tech.prop.2d}. In order to emphasize that we work with a purely 2 dimensional question, we write $(S, \Delta)$  for an  slc, surface  pair. 
It turns out that non-normal pairs such that $\coeff\Delta\subset [\frac12, 1]$ have a rather simple structure.

\begin{thm}\label{slc.cent.12.prop}
 Let $(s\in S, \Delta)$ be a strictly Henselian,  slc, surface  pair 
over an algebraically closed  field of characteristic 0
such that $\coeff\Delta\subset [\frac12, 1]$ and $s\in S$ a non-normal point. 
Then one of the following holds.
\begin{enumerate}
\item The point  $s$ is  an  lc center  and 
  $2(K_S+\Delta)$ is Cartier.  
\item The point  $s$ is  not an  lc center  and $S$ has 
2 irreducible components $(s_i\in S_i, D_i+\Delta_i)$. For both of them 
 the extended dual graph of the minimal  resolution  (over $s\in S$) is of the form
$$
    \xymatrix{%
       \bullet \ar@{-}[r] & c_1 \ar@{-}[r] & \cdots \ar@{-}[r] & c_n  \ar@{-}[r] & \circledast}
    $$
where $\coeff(\bullet)=1$, 
 $\coeff(\circledast)\in [\frac12, 1)$,  $n\geq 0$ and  $c_i\geq 2$ for every $i$. Furthermore  $\diff_{D_1}\Delta_1=\diff_{D_2}\Delta_2$. The local class group has rank 1. 
\item  The point $s$ is not an lc center  and  $S$ has a quasi-\'etale double cover as in {\rm (2)}. The local class group is torsion.
\end{enumerate}
\end{thm}

{\it Note.} The Theorem should hold more generally whenever the characteristic is not 2, but there may be a lack of  references related to adjunction. 
In characteristic 2 there should be only one more case for which the normalization induces an inseparable  map on the conductors. 
\medskip

Proof. 
Let $(s_i\in S_i, D_i+\Delta_i)$ be the  irreducible components of the normalization of $(S, \Delta)$, where the $D_i$ denote the conductors. The pairs $(s_i\in S_i, D_i+\Delta_i)$  are lc and $D_i\neq 0$ for every $i$  since $S$ is not normal at $s$.

The pairs $(s_i\in S_i, D_i+\Delta_i)$ are  described in
(\ref{lc.cent.12.say}.1--3). Correspondingly, there are 2 cases.

{\it Plt case.}  If $s$ is  not an  lc center  then none of the
$s_i$ is an  lc center by \cite[5.10.3]{kk-singbook}.  Thus each 
$(s_i\in S_i, D_i+\Delta_i)$ is as in (\ref{lc.cent.12.say}.1). In particular,
$D_i$ is irreducible and  there are at most 2 irreducible components. 
If there are 
 2  irreducible components  then the gluing is done by an isomorphism
$\tau:D_1\to D_2$.  These give case (2). Any element of the local class group is give by a pair of divisors  $C_i\subset S_i$. Then
$\delta(C_1, C_2):=(C_1\cdot D_1)-(C_2\cdot D_2)$ is well defined and we get an exact sequence
$$
0 \to (\mbox{torsion subgroup})\to \cl(s\in S)\stackrel{\delta}{\longrightarrow} \q.
$$

If  there is  1 irreducible component
 $(s_1\in S_1, D_1+\Delta_1)$ then  the gluing is done by an involution
$\tau:D_1\to D_1$. Furthermore,   it has a quasi-\'etale double cover
$\rho:  (\tilde S, \tilde \Delta)\to (S, \Delta)$,
with 2 irreducible components  by  \cite[5.23]{kk-singbook}, giving case (3). 
The covering involution $\rho$ interchanges the 2 irreducible components
$\tilde S_i$, hence  $\delta\circ \rho=-\delta$. Thus only the
torsion subgroup of $\cl(\tilde s\in \tilde S)$ descends to give divisors on $S$. 

{\it Non-plt case.}  If $s$ is   an  lc center  then  each
$s_i$ is an  lc center by \cite[5.10.3]{kk-singbook}. 
  Thus the irreducible components $(S_i, D_i+\Delta_i)$ are as in
(\ref{lc.cent.12.say}.2--3), their number can be arbitrary. 
Let $D_{ij}$ denote  the irreducible components  of the normalization of  $D_i$.
Thus we have $\rho_i:D_{i1}\amalg D_{i2}\to D_i$ in case (\ref{lc.cent.12.say}.2)
and $\rho_i:D_{i1}\cong D_i$ in case (\ref{lc.cent.12.say}.3).

Let $D:=\amalg_{ij} (s_{ij}\in D_{ij})$ be their disjoint union and let $\tau:D\to D$  denote  the gluing involuton $\tau$ acting on it. Furthermore, by (\ref{lc.cent.12.say}.4), each   $\omega_{S_i}^{[2]}(2D_i+2\Delta_i)$ is a line bundle on $S_i$ and  the Poicar\'e residue map gives canonical isomorphisms
$$
\res_{S_i/D_{ij}}^2: \rho_{ij}^*\omega_{S_i}^{[2]}(2\diff_{D_{ij}}\Delta_i)\cong 
\omega_{D_{ij}}^2\bigl(2\diff_{D_{ij}}\Delta\bigr)
\cong \omega_{\bar D_{ij}}^2\bigl(2[s_{ij}]).
$$
Applying the Poicar\'e residue map twice gives  canonical isomorphisms
$$
\res_{S_i/s_{ij}}^2:\omega_{S_i}^{[2]}(2D_i+2\Delta_i)|_{s_{ij}}\cong k(s_{ij}).
$$
We can thus pick $\tau$-invariant sections 
$$
(\sigma^D_{ij})\in \oplus_{ij} H^0\bigl( D_{ij}, \omega_{D_{ij}}^2\bigl(2\diff_{D_{ij}}\Delta_i\bigr)\bigr)
$$
that have residue 1 at the points $s_{ij}$. Since they have the same residue, they descend to sections
$$
(\sigma^D_{i})\in \oplus_{i} H^0\bigl( D_{i}, \omega_{D_{i}}^2\bigl(2\diff_{D_{i}}\Delta_i\bigr)\bigr).
$$
We can lift these $\sigma^D_{i} $ back to
sections
$$
(\sigma_i)\in \oplus_i H^0\bigl(S_i, \omega_{S_i}^{[2]}(2D_i+2\Delta_i)\bigr).
$$
By \cite[5.8]{kk-singbook} the $(\sigma_i) $ descend to a 
section of $\omega_{S}^{[2]}(2\Delta)$ that has residue 1 (hence nonzero) at the origin. Therefore   $\omega_{S}^{[2]}(2\Delta)$  is locally free at $p$. This completes the proof of (1). \qed

\begin{reminder}\label{lc.cent.12.say}
The list of all lc surface pairs  $(S,\Theta)$ where
$\coeff\Theta\subset [\frac12, 1]$ is given in
\cite[pp.125-128]{kk-singbook}. Here we are interested in those special cases 
when  there is a divisor with coefficient 1. Since the list in \cite{kk-singbook} is organized differently, 
the classification is summarised next where we use  $\Theta:=D+\Delta$. 

 Thus let $(S, \Theta)$ be an lc  pair of dimension 2 over an algebraically closed field such that $\coeff\Theta\subset [\frac12, 1]$ and $s\in \rdown{\Theta}$ a  closed point.
We discribe $(s\in S, \Theta)$ using 
 the extended dual graph of the minimal embedded resolution   $\pi:  S'\to S$.
The exceptional divisors are denoted by the negative of their self-interscetion; $2$ or $c_i$ in the diagrams. 
The birational transforms of the local branches of 
 $\Theta$  are denoted by  $\bullet$ if $\coeff(\bullet)=1$ and
by $\circledast$ if $\coeff(\circledast)\in [\frac12, 1)$. 

There are 3 distinct cases.

{\it Plt case.} If $s$ 
 is not an lc center then the extended dual graph is one of the following, where   $n\geq 0$ and  $c_i\geq 2$ for every $i$.
$$
    \xymatrix{%
       \bullet \ar@{-}[r] & c_1 \ar@{-}[r] & \cdots \ar@{-}[r] & c_n  \ar@{-}[r] & \circledast}
      \eqno{(\ref{lc.cent.12.say}.1)}
    $$

{\it Cyclic non-plt case.} Here $s$ is an lc center,  $n\geq 0$,  $c_i\geq 2$ and the extended dual graph is
$$
    \xymatrix{%
       \bullet \ar@{-}[r] & c_1 \ar@{-}[r] & \cdots \ar@{-}[r] & c_n  \ar@{-}[r] & \bullet}
      \eqno{(\ref{lc.cent.12.say}.2)}
    $$

{\it Dihedral non-plt case.} Here $s$ is an lc center,
 $\coeff(\circledast)=\frac12$, $n\geq 0$ and  $c_i\geq 2$ except that $c_n=1$ is allowed in cases (\ref{lc.cent.12.say}.3.2--3)
$$
    \xymatrix{%
      && & 2 & \\
      \bullet \ar@{-}[r] & c_1 \ar@{-}[r] & \cdots \ar@{-}[r] & c_n \ar@{-}[u]
       \ar@{-}[r] & 2  }
       \eqno{(\ref{lc.cent.12.say}.3.1)}
     $$
$$
    \xymatrix{%
      && & \circledast & \\
      \bullet \ar@{-}[r] & c_1 \ar@{-}[r] & \cdots \ar@{-}[r] & c_n \ar@{-}[u]
       \ar@{-}[r] & 2  }
       \eqno{(\ref{lc.cent.12.say}.3.2)}
     $$
$$
    \xymatrix{%
      && & \circledast & \\
      \bullet \ar@{-}[r] & c_1 \ar@{-}[r] & \cdots \ar@{-}[r] & c_n \ar@{-}[u]
       \ar@{-}[r] & \circledast  }
       \eqno{(\ref{lc.cent.12.say}.3.3)}
     $$
Next let $\Sigma$ 
be the divisor on $S'$ that contains the curves marked by  $\bullet$ or  $c_i$ with coefficient 1 and the curves marked by  $2$ or $\circledast$ with coefficient $\frac12$. 
By inspection  we see that $2(K_{S'}+\Sigma)$ is a $\z$-divisor that has 0 intersection with all $\pi$-exceptional divisors. Thus  $2(K_S+\Theta)=2\cdot \pi_*(K_{S'}+\Sigma)$ is Cartier near $s$ by
\cite[10.9.2]{kk-singbook}. 
\medskip

{\it Conclusion \ref{lc.cent.12.say}.4.}  If $p$ is an lc center then $2(K_S+\Theta)$ is Cartier near $s$. \qed  
\end{reminder}

\begin{exmp}
The nice dichotomy of Theorem~\ref{slc.cent.12.prop} does not seem to
carry over to higher dimensions. As an example, pick
$\frac12\leq c_1, c_2, c_3\leq 1$ such that $c_1+c_2+c_3=2$. 
Consider the pair
$$
\bigl((x_3x_4=0), c_1(x_1=0)+c_2(x_2=0)+c_3(x_1=x_2)\bigr)\subset \a^4.
$$
It is  non-normal  and the origin is an lc center.
\end{exmp}

\section{The Poincar\'e residue map}\label{sec.4}

We study the surjectivity of the  Poincar\'e residue map for slc surface pairs. First we show surjectivity for the pairs listed in (\ref{lc.cent.12.say}.1).
We stress that this is a rather special property of such pairs. 
We see in Example~\ref{diff.nonr.exmp.1} that it fails for some dihedral pairs, even when $\Delta=0$. 
Also, even on smooth surfaces, it fails for every other $\Delta'$ for some $m$; see  Example~\ref{key.plt.exmp.noneq}.

\begin{prop}\label{key.plt.exmp}
Let $( S, D+\Delta)$ be an lc surface pair as in {\rm (\ref{lc.cent.12.say}.1)}, over a field of characteristic 0.
Then  the Poincar\'e residue map
$$
\res_{S/D}^m: 
\omega_S^{[m]}\bigl(mD+\rdown{m\Delta}\bigr)\to
\omega_D^{[m]}\bigl(\rdown{m\diff_D\Delta}\bigr)
\eqno{(\ref{key.plt.exmp}.1)}
$$
is surjective for every $m$.
\end{prop}

Proof. We use the (\'etale-local) respresentation of $( S, D+\Delta)$ as a quotient
$$
\bigl( S, D+\Delta\bigr):=\bigl( \tilde S, \tilde D+\tilde \Delta\bigr)/\tfrac1{n}(1,q),
\eqno{(\ref{key.plt.exmp}.2)}
$$
where $\bigl( \tilde S, \tilde D+\tilde \Delta\bigr):=\bigl(\a^2_{xy}, (y=0)+(1-c)(x=0)\bigr)$; cf.\ \cite[3.32]{kk-singbook}.

We can write the sections of $\omega_S^{[m]}\bigl(mD+\rdown{m\Delta}\bigr) $  in the form
$$
g(x,y)\bigl(\tfrac{dx}{x}\wedge \tfrac{dy}{y}\bigr)^{\otimes m},
\eqno{(\ref{key.plt.exmp}.3)}
$$
where $x^{mc}$ divides $g(x,y)$  (in the ring of Puiseux series)
and $g(x, y)$ is $\mu_n$-invariant. 
Any monomial in $g$ that contains $y$ restricts to 0 on $D$,
thus in (\ref{key.plt.exmp}.2) only the sections of the form
$$
x^r\bigl(\tfrac{dx}{x}\wedge \tfrac{dy}{y}\bigr)^{\otimes m}
$$
have non-zero image.
We need the $\mu_n$-invariant generator, which is
$$
\sigma:=x^{n\rup{mc/n}}\bigl(\tfrac{dx}{x}\wedge \tfrac{dy}{y}\bigr)^{\otimes m}.
\eqno{(\ref{key.plt.exmp}.4)}
$$
Setting $\gamma:=\tfrac{c}{n}$, the sign convention of
\cite[4.1]{kk-singbook} gives that 
$$
\res_{\tilde S/\tilde D}^m(\sigma)=(-1)^mx^{n\rup{m\gamma}}\bigl(\tfrac{dx}{x}\bigr)^{\otimes m}.
\eqno{(\ref{key.plt.exmp}.5)}
$$
On $D$ the local coordinate is $z=x^n$ and
$\tfrac{dx}{x}=n \tfrac{dz}{z} $, hence we get that
$$
\res_{S/D}^m(\sigma)
=
(-n)^mz^{\rup{m\gamma}}\bigl(\tfrac{dz}{z}\bigr)^{\otimes m}.
\eqno{(\ref{key.plt.exmp}.6)}
$$
The different is computed by the formula (cf.\ \cite[3.45]{kk-singbook})
$$
\diff_D\Delta=\bigl(1-\tfrac1{n}+\tfrac{1-c}{n}\bigr)[s]=
\bigl(1-\tfrac{c}{n}\bigr)[s]=(1-\gamma)[s], 
\eqno{(\ref{key.plt.exmp}.7)}
$$
where $s\in D\subset S$ denotes the origin.
Since $m-\rup{m\gamma}=\rdown{m(1-\gamma)}$, (\ref{key.plt.exmp}.6) shows the isomorphism (modulo torsion supported at  $s$)
$$
\res_{S/D}^m: 
\omega_S^{[m]}\bigl(mD+\rdown{m\Delta}\bigr)|_D\congt
\omega_D^{[m]}\bigl(\rdown{m\diff_D\Delta}\bigr). \qed
\eqno{(\ref{key.plt.exmp}.8)}
$$

Next we compute  the $\q$-divisor version of (\ref{key.plt.exmp}.1).

\begin{exmp} \label{pre.key.plt.exmp}
Consider 2 surface pairs
$$
\bigl( S_i, D_i+(1-c_i)C_i\bigr):=\bigl(\a^2_{xy}, (y=0)+(1-c)(x=0)\bigr)/\tfrac1{n_i}(1,1)
\eqno{(\ref{pre.key.plt.exmp}.1)}
$$
and glue them using an isomorphism $\tau:D_1\to D_2$ to get
$$(S,  \Delta):=\bigl( S_1, D_1+(1-c_1)C_1\bigr)\amalg_{\tau}\bigl( S_2, D_2+(1-c_2)C_2\bigr).
\eqno{(\ref{pre.key.plt.exmp}.2)}
$$
Note that $\diff_{D_i}(1-c_i)C_i=1-\tfrac{c_i}{n_i}$, thus
$(S,  \Delta)$ is slc $\Leftrightarrow$ $K_S+\Delta$ is $\q$-Cartier 
$\Leftrightarrow$   $\tfrac{c_1}{n_1}=\tfrac{c_2}{n_2}$.

Given any $n_1, n_2$, choose the $c_i$ such that $\tfrac{c_1}{n_1}=\tfrac{c_2}{n_2}$ and  $c_i<\frac12$. Then 
$2K_S+\rdown{2\Delta}=2K_S+C_1+C_2$. 
Note that $\bigl(2K_S+C_1+C_2\bigr)|_{S_i}=2K_{S_i}+2D_i+C_i$ and
$$
\bigl(2K_{S_i}+2D_i+C_i\bigr)|_{D_i}=2\bigl(K_{D_i}+\bigl(1-\tfrac1{n_i}\bigr)[s]\bigr)+\tfrac1{n_i}[s]=2K_{D_i}+\bigl(2-\tfrac1{n_i}\bigr)[s],
$$
where $s\in D_i\in S_i$ denotes the origin. 
Thus $2K_S+\rdown{2\Delta}$ is  $\q$-Cartier  iff $n_1=n_2$. 
\end{exmp}

Formula (\ref{key.plt.exmp}.8) and
Example~\ref{pre.key.plt.exmp}   directly imply the following.

\begin{cor} \label{key.plt.exmp.cor}
Let $(S,\Delta)=(S_1,D_1+\Delta_1)\amalg_{\tau} (S_2,D_2+\Delta_2)$ be an slc surface as in (\ref{slc.cent.12.prop}.2).  Then
\begin{enumerate}
\item 
$\omega_{S_1}^{[m]}\bigl(mD_1+\rdown{m\Delta_1}\bigr)|_{D_1}\congt
\omega_{S_2}^{[m]}\bigl(mD_2+\rdown{m\Delta_2}\bigr)|_{D_2}$, but in general
\item $ \bigl(mK_{S_1}+mD_1+\rdown{m\Delta_1}\bigr)|_{D_1} \neq
\bigl(mK_{S_2}+mD_2+\rdown{m\Delta_2}\bigr)|_{D_2}$.  \qed
\end{enumerate}
\end{cor}

Next we compute the Poincar\'e residue map in the dihedral cases. 
Although these are not needed for the proof of Theorem~\ref{half.pgsh.2d.inv.thm}, they show that the
Poincar\'e residue map is not surjecrive in general. This suggests that it may not be easy to understand the pluricanonical sheaves  for non-normal pairs using the normalization.  

\begin{exmp} \label{diff.nonr.exmp.1}
Let $(S, B)$ be a pair as in the  dihedral case (\ref{lc.cent.12.say}.2.1). It can also be   obtained as the  quotient
of the pair
$$
\bigl(\tilde S, \tilde B\bigr):=\bigl((xy=z^{2n}), (z=0)\bigr)
\eqno{(\ref{diff.nonr.exmp.1}.1)}
$$
by the involution
$\tau: (x,y,z)\mapsto (y,x,-z)$. Both the $x$ and $y$  axes map isomorphically to $B\subset S$.

A local generator of $\omega_{\tilde S}$ is $z^{-2n+1}dx\wedge dy$, thus 
a local generator of $\omega_{\tilde S}(\tilde B)$ is $z^{-2n}dx\wedge dy$, which can be rewritten as
$$
\sigma:=\tfrac{dx}{x}\wedge \tfrac{dy}{y}.
$$
Thus we see that $\tau^*\sigma=-\sigma$ and so $\sigma$ does not descend to
$S$. Hence $\omega_S(B)$ is not locally free. However
 $\tau^*\sigma^{\otimes 2}=\sigma^{\otimes 2}$, so $\sigma^{\otimes 2}$ does  descend to a local generator of $\omega_S^{[2]}(2B)$, which is thus  locally free.

We can also find generators of $\omega_S(B)$, obtained from the
$\tau$-invariant forms
$$
(x-y)\tfrac{dx}{x}\wedge \tfrac{dy}{y}\qtq{and}
z\tfrac{dx}{x}\wedge \tfrac{dy}{y}.
\eqno{(\ref{diff.nonr.exmp.1}.2)}
$$
By restricting these first to the $x$-axis and then using that the latter is  isomorphic to $B$, we see that 
$$
\res_{S/B}\bigl( (x-y)\tfrac{dx}{x}\wedge \tfrac{dy}{y}\bigr)=-dx
\qtq{and} \res_{S/B}\bigl( z\tfrac{dx}{x}\wedge \tfrac{dy}{y}\bigr)=0.
\eqno{(\ref{diff.nonr.exmp.1}.3)}
$$
Thus $\res_{S/B}:\omega_S(B)|_B\to\omega_B$
 is an isomorphism modulo torsion supported at the origin $s\in S$. 
By contrast, in degree 2 we have  the Poincar\'e residue isomorphism 
$$
\res_{S/B}^2: \omega_S^{[2]}(2B)|_B\cong \omega_B^{2}(2[s]).
\eqno{(\ref{diff.nonr.exmp.1}.4)}
$$
This also shows that
$\diff_B (K_S+B) =K_B+ [s]$.  Thus  the Poincar\'e residue map in degree $m$ gives surjections
$$
\begin{array}{ccccll}
\res_{S/B}^m&:& \omega_S^{[m]}(mB) & \onto & \bigl(\omega_B([s])\bigr)^{\otimes m} & \qtq{if $m$ is even, but}\\
\res_{S/B}^m&:& \omega_S^{[m]}(mB) & \onto & \bigl(\omega_B([s])\bigr)^{\otimes m}(-[s]) & \qtq{if $m$ is odd.}
\end{array}
\eqno{(\ref{diff.nonr.exmp.1}.5)}
$$
\end{exmp}

\begin{exmp} \label{diff.nonr.exmp.2}
Consider 2 pairs  $(S_1, B_1)$  and  $(S_2, B_2)$  as in the  dihedral case (\ref{lc.cent.12.say}.2.1) and
$(T, C_1+C_2):=(\a^2_{\mathbf u}, (u_1=0)+(u_2=0)\bigr)$. 
Out of these we can assemble a reducible slc pair
by gluing $S_i$ to $T$ using  isomorphisms $\tau_i:C_i\to B_i$ such that $\tau_i^*dx_i=du_i$. We get an slc pair
$$
(S, 0)\cong (S_1, B_1)\amalg_{\tau_1}(T, C_1+C_2)\amalg_{\tau_2}(S_2, B_2).
\eqno{(\ref{diff.nonr.exmp.2}.1)}
$$
We claim that, modulo torsion supported at the origin $s\in S$,
\begin{enumerate}\setcounter{enumi}{1}
\item $\omega_S|_{S_i}\congt \omega_{S_i}(B_i)$ but
\item $\omega_S|_{T}\congt \omega_{T}(C_1+C_2)(-[s])$.
\end{enumerate}
Thus, although $\omega_S$ is $S_2$, its restriction to
$T$ is not $S_2$. This makes it hard to study the depth of pluricanonical sheaves on reducible slc pairs using the normalization.

In order to see (2) we need to show that every section of 
$ \omega_{S_i}(B_i)$ extends to a section of $\omega_S$.
It is enough to prove this for the generators in
(\ref{diff.nonr.exmp.1}.2). 
Here $ z_i\tfrac{dx_i}{x_i}\wedge \tfrac{dy_i}{y_i}$ vanishes on $B_i$, so we can extend it by 0 to the other components. The other  generator
$(x_i-y_i)\tfrac{dx_i}{x_i}\wedge \tfrac{dy_i}{y_i}$
restricts to $B_i$ as $-dx_i=-du_i$. This can be extended to $T$ as
$u_i\tfrac{du_i}{u_i}\wedge \tfrac{du_{3-i}}{u_{3-i}}$ and then as zero to $S_{3-i}$. 

There is a sign ambiguity in the definition of higher codimension restriction maps, but
 $\pm \res_{S/C_i}=\pm \res_{S/B_i}=\pm \res_{S_i/B_i}\circ \res_{S/S_i}$, hence 
$\res_{S/C_i}$ gives  a surjection
$\res_{S/C_i}:\omega_S\onto \omega_{C_i}$. 
On the other hand, $\res_{T/C_i}$ gives  a surjection
$\res_{T/C_i}: \omega_{T}(C_1+C_2)\onto \omega_{C_i}([s])$. These 
 show that the image of 
$\omega_S|_{T}$ is contained in $ \omega_{T}(C_1+C_2)(-[s])$.
The latter is generated by the forms $u_i\tfrac{du_i}{u_i}\wedge \tfrac{du_{3-i}}{u_{3-i}}$  and 
we already saw that these extend to sections of $\omega_S$. This proves (3).

\end{exmp}

The next example shows that Proposition~\ref{key.plt.exmp}
does not hold if $\Delta$ has at least 2 irreducible components.

\begin{exmp} \label{key.plt.exmp.noneq}
Let $D, C_1, \dots, C_r$ be distinct lines through the origin $s\in S:=\a^2$. For some positive rational numbers $c_i$ consider the pair
$(S, D+\Delta)$ where $\Delta=\sum_i c_iC_i$.
We claim that if $r\geq 2$ then there is an $m>0$ such that
$$
\res_{S/D}^m: \omega_S^m(\rdown{m\Delta})|_D\to  \omega_D^m(\rdown{m\diff_D\Delta})
\eqno{(\ref{key.plt.exmp.noneq}.1)}
$$
is not surjective.

First note that $\diff_D\Delta=(\sum_i c_i)[s]$ and
the 2 sides of (\ref{key.plt.exmp.noneq}.1) are
$$
\res_{S/D}^m: \omega_D^m\bigl(\tsum_i \rdown{mc_i}\cdot [s]\bigr)\to 
\omega_D^m\bigl(\rdown{m\tsum_i c_i}\cdot [s]\bigr).
\eqno{(\ref{key.plt.exmp.noneq}.2)}
$$
Thus our claim is equivalent to saying that
$$
\tsum_i \rdown{mc_i} < \rdown{m\tsum_i c_i} \qtq{for some} m>0.
\eqno{(\ref{key.plt.exmp.noneq}.3)}
$$
First choose the smallest $m>0$ such that $m\sum_i c_i$ is an integer.  Then $\sum_i \rdown{mc_i}\leq  \sum_i mc_i=\rdown{m\tsum_i c_i}$ and equality holds iff all the $mc_i$ are integers. Thus we are done unless $\sum_i c_i=a/m$  for some $(a,m)=1$
and $c_i=a_i/m$ for some integers $a_i$. 

Now choose $m'>0$ such that $m'a\equiv 1\mod m$. Note that $m'a_i/m$ is not an integer since $(m',m)=1$ and $c_i<1$.  Thus 
$$
\tsum_i \rdown{m'c_i}\leq \sum_i \bigl(m'c_i-\tfrac1{m}\bigr)
= m'(\tsum_i c_i)-\frac{r}{m}< m'(\tsum_i c_i)-\frac1{m}=\rdown{m'\tsum_i c_i}.
$$

\end{exmp}

\section{Standard coefficients}\label{sec.5}

In this section we prove Corollary~\ref{st.pg.inv.cor}. More generally, 
 we study what happens in all dimensions if $\coeff\Delta$ is contained in the
{\it standard coefficient set}  ${\mathbf T}:=\{\frac12, \frac23, \frac34, \dots,1\}$.
The first observation is that   then a strong form of  Proposition~\ref{tech.prop.2d} holds in all dimensions. 
We write $(X,\Delta)$ for pairs of arbitrary dimension but change to
$(S,\Delta)$ when the discussion is restricted to surfaces.

\begin{prop} \label{tech.prop.std.prop}
Let $(X,\Delta)$ be an  slc pair  over a field of characteristic 0 where $\coeff\Delta\subset \{\frac12, \frac23, \frac34, \dots,1\}$. 
Assume that $x\in X$ is not an lc center. 
Then 
$$
\depth_x\omega_X^{[m]}\bigl(\rdown{m\Delta}-B\bigr)\geq \min\{3, \codim_Xx\}.
\eqno{(\ref{tech.prop.std.prop}.1)}
$$ 
for every $m\in \z$ and $B\subset \rdown{\Delta}$.
\end{prop}

Proof. Set $-D\sim  mK_X+\rdown{m\Delta}-B$ and note that
$$
D\simq -m(K_X+\Delta)+\{m\Delta\}+B\simq  \{m\Delta\}+B\leq \Delta.
$$
Thus (\ref{tech.prop.std.prop}.1) follows from \cite{k-dep}; see also \cite[1.81]{kk-singbook}. \qed

\begin{cor} \label{tech.prop.std.cor}
Let $(X,H+\Delta)$ be an  slc pair  over a field of characteristic 0 where $H$ is Cartier and $\coeff\Delta\subset \{\frac12, \frac23, \frac34, \dots,1\}$. 
Then 
$$
\depth_x\omega_X^{[m]}\bigl(\rdown{m\Delta}-B\bigr)\geq \min\{3, \codim_Xx\}.
$$
holds for every $m\in \z$,  $B\subset \rdown{\Delta}$ and $x\in H$.
\end{cor}

Proof.  By the monotonicity of discrepancies  \cite[2.27]{km-book}, none of the lc centers of $(X,\Delta)$  is contained in $H$. Thus the Corollary follows from  
Proposition~\ref{tech.prop.std.prop}. \qed
\medskip

As we noted at the beginning of Section~\ref{sec.1}, Corollary~\ref{tech.prop.std.cor} and 
 \cite[Prop.16]{k-lpg1} imply the following.

\begin{cor}  \label{st.pg.inv.cor.ndim}
Let $S$ be a reduced scheme over a field of characteristic 0 and $f:(X,\Delta) \to S$ a stable morphism   such that $\coeff\Delta\subset \{\frac12, \frac23, \frac34, \dots, 1\}$. Then, for every $m\geq 1$,
$\omega_{X/S}^{[m]}\bigl(\rdown{m\Delta}\bigr)$
is flat over $S$ and commutes with base change. \qed
\end{cor}

By the Cohomology and Base Change Theorem, the $n$-dimensional version of
Corollary~\ref{st.pg.inv.cor} follows once we establish vanishing theorems for the fibers $\omega_{X_s}^{[m]}\bigl(\rdown{m\Delta_s}\bigr) $.
More generally, the following should be true.

\begin{conj} \label{stand.pluri.vanish.conj}
Let $(X,\Delta)$ be an proper,  slc pair  over a field of characteristic 0 such that $K_X+\Delta$ is ample.  Fix $m\geq 2$ and assume that $$
\coeff\Delta\subset \{\tfrac12, \tfrac23, \tfrac34, \dots,1\}\cup \bigl[1-\tfrac1m, 1\bigr].
\eqno{(\ref{stand.pluri.vanish.conj}.1)}
$$ 
Then 
$$
H^i\bigl(X, \omega_X^{[m]}(\rdown{m\Delta})\bigr)=0\qtq{for}  i>0.
\eqno{(\ref{stand.pluri.vanish.conj}.2)}
$$
\end{conj}

Proof attempt.  Note that
$$
mK_X+\rdown{m\Delta}\simq K_X+(m-1)(K_X+\Delta)+
\rdown{m\Delta}-(m-1)\Delta.
\eqno{(\ref{stand.pluri.vanish.conj}.3)}
$$
Our assumption (\ref{stand.pluri.vanish.conj}.1) guarantees that
$0\leq \rdown{m\Delta}-(m-1)\Delta\leq \Delta$. 
If $mK_X+\rdown{m\Delta}$ is $\r$-Cartier then we can apply 
the Ambro-Fujino form of Kodaira's vanishing 
(see Theorem~\ref{fujino.1.10.thm})  and we are done. However, usually 
$mK_X+\rdown{m\Delta}$ is not $\r$-Cartier. 
It is a natural idea to try to find a proper, birational  morphism
$\pi: (X', \Delta')\to (X, \Delta)$ such that
 $mK_{X'}+\rdown{m\Delta'}$  is $\r$-Cartier, 
establish vanishing on $X'$ and then  descend to $X$.
That is, we aim to find a proper, birational  morphism
 $\pi: (X', \Delta')\to (X, \Delta)$ 
 with the following properties.
\begin{enumerate}\setcounter{enumi}{3}
\item $\pi_*\bigl(\omega_{X'}^{[m]}(\rdown{m\Delta'})\bigr)=
\omega_{X}^{[m]}(\rdown{m\Delta})$, 
\item $H^i\bigl(X', \omega_{X'}^{[m]}(\rdown{m\Delta'})\bigr)=0$ for $i>0$ and
\item $R^i\pi_*\bigl(\omega_{X'}^{[m]}(\rdown{m\Delta'})\bigr)=0$ for $i>0$.
\end{enumerate}
Then  the Leray spectral sequence  shows that
$H^i\bigl(X, \omega_{X}^{[m]}(\rdown{m\Delta})\bigr)=0$ for $i>0$. 

If  $X$ is normal, then \cite[Prop.19]{k-lpg1} and Theorem~\ref{fujino.1.10.thm}
show that there is a small
modification $\pi:X'\to X$ with these properties. 
However, if $X$ is not normal, then sometimes there is no such small
modification. 
This is obvious for surfaces, since 
a demi-normal surface has no nontrivial small modifications. Therefore  we have to use a  birational  morphism with exceptional divisors. This brings in 2 extra problems. 
\begin{itemize}
\item  In many cases the coefficient of an exceptional divisor  in $\Delta'$ should be the discrepancy, or a small perturbation of it. Thus it  may not satisfy the numerical assumptions  (\ref{stand.pluri.vanish.conj}.1).
\item If an exceptional divisor appears with  coefficient 1, then the needed  vanishing claims (\ref{stand.pluri.vanish.conj}.5--6) usually do  not hold.
\end{itemize} 
For surfaces we can avoid  these problems, but only with  very special choices of $\pi$.

\begin{say}[Proof of Conjecture~\ref{stand.pluri.vanish.conj}  for surfaces]
\label{pf.of.stand.pluri.vanish.conj}
There are only finitely many points $s\in S$ such that
$mK_{S}+\rdown{m\Delta}$ is  not $\r$-Cartier at $s$. At these points $S$ is non-normal. We follow the classification of such points given in
(\ref{slc.cent.12.prop}.1--3) and in each case give a local description of $\pi:(S', \Delta')\to (S, \Delta)$. We describe $\pi$ after passing to the strict Henselisation, but in each case this automatically descends to the original base field.

\medskip
{\it \ref{pf.of.stand.pluri.vanish.conj}.1 (Plt case with 2 components.)}  By  (\ref{slc.cent.12.prop}.2)   here $(S, \Delta)$ is glued together from
2 branches  $S_1, S_2$ , with resolution dual graphs
$$
    \xymatrix{%
       \bullet \ar@{-}[r] & c_{i1} \ar@{-}[r] & \cdots \ar@{-}[r] & c_{in_i}  \ar@{-}[r] & \circledast}
    $$
where $\coeff(\bullet)=1$, 
 $\coeff(\circledast)=1-d_i$. We choose the partial resolution
$S'_i\to S_i$ that  extracts only the curve  $C_{i1}$. Thus $S'_i$ has 1 singular point
(obtained by contracting the curves $C_{i2},\dots, C_{in_i}$)  and
$D_i\cap C_{i1}$ is a smooth point of $S'_i$. 
Note further that 
$$
a(C_{i1}, S_i, \Delta_i)=-1+\tfrac1{n_i}-\tfrac{1-d_i}{n_i}=-1+\tfrac{d_i}{n_i}.
$$
As we noted in Example~\ref{pre.key.plt.exmp}, the quatity  $\gamma:=\tfrac{d_i}{n_i}$ is independent of $i$.
We add $C_{i1}$ to  $\Delta'_i$  with coefficient  $1-\gamma$,
thus $\Delta'_i$ aso satisfies the coefficient assumption
(\ref{stand.pluri.vanish.conj}.1).
 
We can now glue $S'_1$ and $S'_2$ to get $S'\to S$.
On $S'$ we get 2 normal cyclic quotient singularities and 1 normal crossing point  of the form 
$\bigl((xy=0), (1-\gamma)(z=0)\bigr)$. 
\medskip

{\it  \ref{pf.of.stand.pluri.vanish.conj}.2 (Plt case with 1 component.)}  By (\ref{slc.cent.12.prop}.3) in these cases the local class group is torsion, so $mK_{S}+\rdown{m\Delta}$ is   $\r$-Cartier at $s$. 
\medskip

{\it  \ref{pf.of.stand.pluri.vanish.conj}.3 Log center case.}  
If  $m\Delta$ is a $\z$-divisor  then $mK_S+\rdown{m\Delta}=mK_S+m\Delta$ is $\q$-Cartier by (\ref{slc.cent.12.prop}.1). 
Thus assume that  $\rdown{m\Delta}\lneqq m\Delta$. 
For $\pi: S'\to S$ we take the slc modification. That is, 
on the irreducible components listed in  (\ref{lc.cent.12.say}.3.1--3) we extract all curves marked $c_i$, but we do not extract the curves marked $2$
(these have discrepancy $-\frac12$). Thus all the exceptional curves $C_{ij}$ appear with coefficient 1.  Set $\Delta'':=\pi^{-1}_*\Delta+\sum C_{ij}$.
As we noted, the problem is that we can not apply
Theorem~\ref{fujino.1.10.thm} to  $(S',  \Delta'')$. 

Let $\sigma$ be a section of $\omega_{S}^{[m]}(\rdown{m\Delta})$. We can then view
$\pi^{-1}_*\sigma$ as a section of  $\bigl(\omega_{S'}^{[m]}(\rdown{m\Delta''})\bigr)$. 
Since  $\rdown{m\Delta}\lneqq m\Delta$, the section $\pi^{-1}_*\sigma$
vanishes along all  the exceptional curves $C_{ij}$. Thus we can 
decrease the coefficents of the $C_{ij}$ without violating
(\ref{stand.pluri.vanish.conj}.4). 

To do this choose a $\pi$-exceptional, $\q$-Cartier  divisor $E$ such that
$-E$ is $\pi$-ample  and set
$\Delta':=\Delta''-\epsilon E$.
Then $K_{S'}+  \Delta'$ is  $\pi$-ample on $S'$ for every $\epsilon$.
Furthermore, once we patch the local modifications to a global
$S'\to S$, the divisor  $K_{S'}+  \Delta'$ is  still nef and log big  on $S'$
for  $0<\epsilon\ll 1$. (It has degree 0 only on the $\pi$-exceptional curves over the plt points.)
\medskip

With these choices $\pi: (S', \Delta')\to (S, \Delta)$ 
satisfies (\ref{stand.pluri.vanish.conj}.4) and
(\ref{stand.pluri.vanish.conj}.5--6) follow from
Theorem~\ref{fujino.1.10.thm}. Thus 
 the Leray spectral sequence  shows that
$H^i\bigl(S, \omega_{S}^{[m]}(\rdown{m\Delta})\bigr)=0$ for $i>0$. 
\qed
\end{say}

The following is proved in \cite{ambro} and \cite[1.10]{MR3238112}, see also 
\cite{fuj-book}, where it is called a Reid-Fukuda--type vanishing theorem.
The 2 dimensional case that we use is much easier. 

\begin{thm}[Ambro-Fujino vanishing theorem]
 \label{fujino.1.10.thm}   Let $(X, \Delta)$ be an slc pair and $D$ a Mumford $\z$-divisor on $X$ (that is, $X$ is regular at all generic points of $\supp \Delta$).  Let $f:X\to S$ be a proper morphism. 
Assume that $D\simr K_X+L+\Delta$, where $L$ is $\r$-Cartier, $f$-nef and  log $f$-big.  Then
$R^if_*\o_X(D)=0$ for $ i>0$. \qed
\end{thm}


\def\cprime{$'$} \def\cprime{$'$} \def\cprime{$'$} \def\cprime{$'$}
  \def\cprime{$'$} \def\cprime{$'$} \def\dbar{\leavevmode\hbox to
  0pt{\hskip.2ex \accent"16\hss}d} \def\cprime{$'$} \def\cprime{$'$}
  \def\polhk#1{\setbox0=\hbox{#1}{\ooalign{\hidewidth
  \lower1.5ex\hbox{`}\hidewidth\crcr\unhbox0}}} \def\cprime{$'$}
  \def\cprime{$'$} \def\cprime{$'$} \def\cprime{$'$}
  \def\polhk#1{\setbox0=\hbox{#1}{\ooalign{\hidewidth
  \lower1.5ex\hbox{`}\hidewidth\crcr\unhbox0}}} \def\cdprime{$''$}
  \def\cprime{$'$} \def\cprime{$'$} \def\cprime{$'$} \def\cprime{$'$}
\providecommand{\bysame}{\leavevmode\hbox to3em{\hrulefill}\thinspace}
\providecommand{\MR}{\relax\ifhmode\unskip\space\fi MR }
\providecommand{\MRhref}[2]{%
  \href{http://www.ams.org/mathscinet-getitem?mr=#1}{#2}
}
\providecommand{\href}[2]{#2}

\bigskip

\noindent  Princeton University, Princeton NJ 08544-1000

{\begin{verbatim} kollar@math.princeton.edu\end{verbatim}}

\end{document}